\documentclass[oneside,notitlepage,11pt]{article}

\usepackage{amssymb}
\usepackage[leqno]{amsmath}
\usepackage{amsfonts}
\usepackage{amsopn}
\usepackage{amstext}
\usepackage{amsthm}
\usepackage{enumitem}

\usepackage{tikz-cd}

\usepackage[colorlinks]{hyperref}

\usepackage{calrsfs}
\usepackage{fourier-orns}
\usepackage{clock} 

\providecommand{\cal}{\mathcal}
\renewcommand{\Bbb}{\mathbb}
\renewcommand{\frak}{\mathfrak}
\newenvironment{pf}{\begin{proof}}{\end{proof}}

\newcommand{\Aaa}{{\cal{A}}}

\newcommand{\Dee}{{\cal{D}}}
\newcommand{\Ef}{{\cal{F}}}
\newcommand{\Gee}{{\cal{G}}}

\newcommand{\Yu}{{\cal{U}}}

\newcommand{\Zee}{{\Bbb{Z}}}
\newcommand{\Emm}{{\frak{M}}}

\newcommand{\Nat}{{\Bbb{N}}}
\newcommand{\Qyu}{{\Bbb{Q}}}
\newcommand{\Err}{{\Bbb{R}}}

\newcommand{\lam}{{\lambda}}
\newcommand{\al}{\alpha}

\newcommand{\sig}{\sigma}

\renewcommand{\phi}{\varphi}
\renewcommand{\rho}{\varrho}

\newcommand{\rest}{\restriction}

\newcommand{\ntr}{{n\in\omega}}
\newcommand{\Ntr}{n\in{\Bbb{N}}}
\newcommand{\loe}{\leqslant}
\newcommand{\goe}{\geqslant}

\newcommand{\subs}{\subseteq}
\newcommand{\sups}{\supseteq}
\newcommand{\nnempty}{\ne\emptyset}

\newcommand{\ctbls}[1]{\dpower{#1}{\aleph_0}}

\newcommand{\id}[1]{{\operatorname{id}_{#1}}} 
\newcommand{\cf}{\operatorname{cf}}
\newcommand{\cov}{\operatorname{cov}}
\newcommand{\covsig}{\operatorname{cov}_\sigma}

\newcommand{\oraz}{\qquad\text{and}\qquad}

\newcommand{\poset}{{\Bbb{P}}}

\newcommand{\concat}{{}^\smallfrown}

\newcommand{\setof}[2]{\{#1\colon #2\}}

\newcommand{\sett}[2]{\{#1\}_{#2}}
\newcommand{\sn}[1]{\{#1\}} 
\newcommand{\dn}[2]{\{#1,#2\}} 
\newcommand{\pair}[2]{\langle #1, #2 \rangle} 
\newcommand{\map}[3]{#1\colon #2 \to #3} 
\newcommand{\img}[2]{#1[#2]} 

\newcommand{\dpower}[2]{[#1]^{#2}}

\newcommand{\fra}{Fra\"iss\'e}

\providecommand{\nat}{\omega}

\newcommand{\ciag}[1]{{\sett{{#1}_n}{\ntr}}}

\newcommand{\iso}{\approx}

\newcommand{\ciagi}[1]{\sig{#1}}

\newcommand{\cmp}{\circ} 

\newcommand{\her}[1]{#1\!\!\downarrow\;}

\newcommand{\proto}[1]{{\mathbb S_\kappa}}

\newtheorem{tw}{Theorem}[section]
\newtheorem{wn}[tw]{Corollary}
\newtheorem{lm}[tw]{Lemma}
\newtheorem{prop}[tw]{Proposition}
\newtheorem{claim}[tw]{Claim}
\theoremstyle{definition}
\newtheorem{df}[tw]{Definition}
\newtheorem{ex}[tw]{Example}

\newtheorem{question}[tw]{Question}

\theoremstyle{remark}

\newcommand{\sbs}[1]{{#1}^{\hookleftarrow}}

\newcommand{\BM}[1]{\operatorname{BM}\left(#1\right)}
\newcommand{\BMG}[2]{\BM {#1, #2}}
\newcommand{\BMGE}[2]{\BM {#1, \sbs{#2}}}

\newcommand{\Age}[1]{\operatorname{Age}(#1)}

\title{Games with finitely generated structures}
\author{{\sc Adam Krawczyk}\\
	{\small Insititute of Mathematics}\\
	{\small University of Warsaw, Poland}
\and
{\sc Wies{\l}aw Kubi\'s}\footnote{
	Research of the second author supported by GA\v CR grant No. 17-27844S.}\\
	{\small Institute of Mathematics}\\
	{\small Czech Academy of Sciences}\\
	{\small Prague, Czechia}\\
	------ \\
	{\small Institute of Mathematics}\\
	{\small Cardinal Stefan Wyszy\'nski University}\\
	{\small Warsaw, Poland}
}
\date{\today\ \clocktime}

\begin{document}

\maketitle

\begin{abstract}
We study the abstract Banach-Mazur game played with finitely generated structures instead of open sets. We characterize the existence of winning strategies aiming at a single countably generated structure. We also introduce the concept of \emph{weak \fra\ classes}, extending the classical \fra\ theory, revealing its relations to our Banach-Mazur game.

\ 


\ 

\end{abstract}


\section{Introduction}

Infinite mathematical structures are often built inductively as unions of chains (sometimes called towers or nests) of finite structures of the same type, for example, at each stage adding just one element. Such a procedure could also be random, typical examples come from graph theory where one builds an infinite graph starting from a single vertex, at each stage adding one more, choosing randomly the new edges. Another possibility of building a structure from small ``pieces" is using a natural infinite game in which two players alternately extend given finite structures, building an infinite chain. In the end, after infinitely many steps, one can test the union of the resulting chain. A general question is when one of the players has a strategy for obtaining an object with prescribed properties. Ultimately, the property could be ``being isomorphic to a specific model". A more relaxed version is ``being embeddable into a specific model".

The aim of this note is to characterize the existence of winning strategies in these games. It turns out that the crucial point is the weak amalgamation property, necessary and sufficient for the existence of winning strategies.
We prove that if the goal of the game is ``being isomorphic to $M$" then the second player has a winning strategy if and only if $M$ is a special model characterized by some kind of extension property. Otherwise the first player has a winning strategy, that is, the game is determined.
We also show that if the goal of the game is ``being embeddable into $N$" then the second player has a winning strategy if and only if $N$ contains an isomorphic copy of the special model $M$ mentioned before.
In order to make precise statements of our main results, we now introduce our main concepts.

Throughout this note, $\Ef$ will denote a fixed class of finitely generated models (structures) of a fixed countable signature. We assume that $\Ef$ is closed under isomorphisms. $\ciagi \Ef$ will denote the class of all models representable as unions of countable chains in $\Ef$.

We consider the following infinite game for two players \emph{Eve} and \emph{Odd}. Namely, Eve starts by choosing $A_0 \in \Ef$. Odd responds by choosing a bigger structure $A_1 \sups A_0$ in $\Ef$. Eve responds by choosing $A_2 \in \Ef$ containing $A_1$. And so on; the rules for both players remain the same.
After infinitely many moves, we receive a model $A_\infty = \bigcup_{\ntr}A_n$, called the \emph{result} of the play. Now let $\Aaa \subs \ciagi \Ef$ be a class of models. We can say that \emph{Odd wins} if the resulting structure is isomorphic to some model in $\Aaa$.
Otherwise, \emph{Eve wins}.
The game is particularly interesting when $\Aaa = \sn G$ for some concrete model $G$..
We shall denote this game $\BMG{\Ef}{\Aaa}$ and $\BMG{\Ef}{G}$ in case where $\Aaa = \sn G$.

This is in fact an abstract version of the well known Banach-Mazur game~\cite{Telgarsky}, in which open sets are replaced by abstract objects.
In \cite{KubBM} it was shown that Odd has a winning strategy in $\BMG{\Ef}{G}$ whenever $\Ef$ is a \fra\ class (equivalently: $G$ is homogeneous with respect to its finitely generated substructures). The paper \cite{KubBM} contains also examples of non-homogeneous graphs $G$ for which Odd still has a winning strategy.
One needs to admit that our Banach-Mazur game is a particular case of infinite games often considered in model theory, where sometimes the players are denoted by $\exists$ and $\forall$.
For more information, see the monograph of Hodges~\cite{Hodges-games}.

Our goal is to present a characterization of the existence of a winning strategy in $\BMG{\Ef}{G}$, where $G$ is a countable first-order structure. We also develop the theory of limits of \emph{weak} \fra\ classes, where the amalgamation property is replaced by a weaker condition exhibited first by Ivanov~\cite{Ivanov}, then used by Kechris and Rosendal~\cite{KecRos}, and recently studied by Kruckman~\cite{KruckPhD}. We show that Odd has a winning strategy in $\BMG{\Ef}{G}$ if and only if $\Ef$ is a weak \fra\ class and $G$ is its limit.
We note that similar results, using the topological Banach-Mazur game, were recently obtained by Kruckman in his Ph.D. thesis~\cite{KruckPhD}. Our approach is direct, we do not use any topology. No prerequisites, except a very basic knowledge in model theory, are required for understanding our concepts and arguments.

\paragraph{Prehomogeneity and ubiquity in category.}
As it happens, there have been some works addressing the question when a fixed countable structure occurs residually in a suitably defined space of structures (typically: all structures whose universe is the set of natural numbers).
This idea was first studied by Cameron~\cite{Cam}, later explored by Pouzet and Roux~\cite{PouRou},
although it actually goes back to Pabion~\cite{Pab} who invented and characterized the concept of {prehomogeneity}.
A structure $M$ is \emph{prehomogeneous} if for every finitely generated substructure $A$ of $M$ there is a bigger finitely generated substructure $B$ containing $A$ such that every embedding of $A$ into $M$ extends to an automorphism of $M$ as long as it is extendible to an embedding of $B$ into $M$.
It may easily happen that no automorphism of $M$ extends a given embedding of $B$ into $M$.
Pouzet and Roux proved that a countable structure is ubiquitous in category (i.e. its isomorphic type forms a residual set in a suitable space) if and only if it is prehomogeneous (Cameron~\cite{Cam} showed earlier that homogeneous structures are ubiquitous). The theorem of Pouzet and Roux can be easily derived from our results, due to Oxtoby's characterization~\cite{Oxtoby} of winning strategies in the original Banach-Mazur game played with open sets.

A structure $M$ is \emph{pseudo-homogeneous} (we would rather call it \emph{cofinally homogeneous}) if every finitely generated substructure $A$ of $M$ extends to a finitely generated substructure $B$ of $M$ such that every embedding of $B$ into $M$ extends to an automorphism of $M$. Clearly, this is a natural strengthening of prehomogeneity.
Pseudohomogeneity was, according to our knowledge, first studied by Calais~\cite{Cal1, Cal2}, around ten years after \fra's work~\cite{Fraisse}. The first example of a countable prehomogeneous and not cofinally homogeneous structure is contained in Pabion~\cite{Pab}, attributed to Pouzet. None of these works mentions the weak amalgamation property, an essential tools for constructing prehomogeneous structures from ``small pieces". Also, none of the above-mentioned works combines prehomogeneity or the weak amalgamation property with the existence of winning strategies in the abstract Banach-Mazur game played with finitely generated models.

\section{The setup}

Throughout this note $\Ef$ will always denote a class, closed under isomorphisms, consisting of countable finitely generated structures of a fixed first-order language. We will denote by $\ciagi \Ef$ the class of all structures of the form $\bigcup_{\ntr}X_n$, where $\ciag{X}$ is an increasing chain of structures from $\Ef$.
The relation $X \loe Y$ will mean, as usual, that $X$ is a substructure of $Y$.
We define the \emph{hereditary closure} of $\Ef$ by
$$\her \Ef = \setof{X}{(\exists\; Y \in \Ef)\; X \loe Y \text{ and $X$ is finitely generated }}.$$
Note that $\ciagi(\her \Ef)$ may be strictly larger than $\ciagi \Ef$ (see Example~\ref{ExNonUnivWFs} below).
Recall that $\Ef$ is \emph{hereditary} if $\her \Ef = \Ef$.
Recall that $\Ef$ has the \emph{joint embedding property} (JEP) if for every $X,Y \in \Ef$ there is $Z \in \Ef$ such that $X \loe Z$ and $Y \loe Z$.
A model $M \in \ciagi{\Ef}$ is \emph{$\Aaa$-universal} if every $X \in \Aaa$ embeds into $M$. 
The \emph{age} of a model $M$, denoted by $\Age{M}$, is the class of all finitely generated models embeddable into $M$.
We shall say that $M$ is \emph{universal} if it is $\ciagi{\Ef}$-universal.
Note that $\ciagi{\Ef}$ has an $\Ef$-universal model if and only if $\Ef$ has countably many isomorphic types and the joint embedding property (cf. \cite[Thm. 2.1]{Cam}).
Thus, a hereditary class $\Ef$ is of the form $\Age{M}$ if and only if it has both the JEP and countably many isomorphic types (briefly: CMT).

We will consider the game $\BMG{\Ef}{G}$ described in the introduction, where both players are allowed to play with structures from $\Ef$ and $G \in \ciagi \Ef$.
In Section~\ref{SeccjaSesta} we will consider a relaxed version $\BMGE{\Ef}{G}$, where $\sbs{G}$ denotes the class of all models embeddable into $G$. It is formally easier for Odd to win this game, namely, if Odd has a winning strategy in $\BMG{\Ef}{G}$ then the same strategy is winning in $\BMGE{\Ef}{H}$ for every $H \in \ciagi{\Ef}$ containing $G$.
Nevertheless, we shall show in Section~\ref{SeccjaSesta} that both games are equivalent in the sense that they imply the same properties of $\Ef$.
The following facts are rather straightforward; for the Reader's convenience we provide the proofs.

\begin{prop}\label{Pjepctmts}
	Suppose Odd has a winning strategy in $\BMGE{\Ef}{G}$ for some $G \in \ciagi \Ef$.
	Then $\Ef$ has the JEP and CMT.
\end{prop}

\begin{pf}
	Eve can start with an arbitrary $X \in \Ef$. As Odd has a winning strategy, we deduce that $G$ contains copies of all structures in $\Ef$, while on the other hand $G$ can contain only countably many finitely generated substructures, which shows CMT. The JEP also follows, because given $X,Y \loe G = \bigcup_{\ntr}G_n$ with each $G_n \in \Ef$, we can find $m$ such that $X, Y \loe G_m$, as $X$ and $Y$ are finitely generated.
\end{pf}

\begin{prop}\label{PrpUnqns}
	Let $G_0, G_1$ be such that Odd has winning strategies both in $\BMG{\Ef}{G_0}$ and $\BMG{\Ef}{G_1}$. Then $G_0$ is isomorphic to $G_1$.
\end{prop}

\begin{pf}
	Let $\Sigma_i$ be Odd's winning strategy in $\BMG{\Ef}{G_i}$, where $i=0,1$.
	Let us play the game assuming that Odd is using strategy $\Sigma_0$. 
	Eve starts with some randomly chosen $A_{-1} \in \Ef$ and her choice $A_0$ is computed according to strategy $\Sigma_1$. From that point on, Eve is using strategy $\Sigma_1$ applied to sequences of the form
	$$A_{-1} \loe A_0 \loe \dots \loe A_{n-1},$$
	where $n$ is even.
	In this situation both players win, showing that $A_\infty = \bigcup_{\ntr}A_n$ is isomorphic to both $G_0$ and $G_1$.
\end{pf}

In the sequel we shall frequently use the following trivial fact.

\begin{lm}\label{LMliftinisos}
	Let $A, B \in \Ef$, let $\map h A B$ be an isomorphism, and let $B'\in \Ef$ be such that $B' \goe B$.
	Then there exist $A' \goe A$ and an isomorphism $\map {h'}{A'}{B'}$ extending $h$.
\end{lm}

\section{Weak injectivity}

In this section we exhibit the crucial property of a model $G$ equivalent to the existence of Odd's winning strategy in the game $\BMG{\Ef}{G}$.
As before, $\Ef$ is a fixed class of countable finitely generated structures, and $\ciagi \Ef$ is the class of all unions of countable chains of structures from $\Ef$.
Recall that $G \in \ciagi \Ef$ is \emph{$\Ef$-universal} if every structure from $\Ef$ embeds into $G$.

\begin{prop}\label{PropWEPs}
	Let $\Ef$ be as above and let $G \in \ciagi \Ef$.
	The following properties are equivalent:
	\begin{enumerate}[itemsep=0pt]
		\item[{\rm(a)}] For every $A \loe G$, $A \in \Ef$, there is $B \in \Ef$ such that $A \loe B \loe G$ and for every $X \in \Ef$ with $B \loe X$ there exists an embedding $\map{f}{X}{G}$ satisfying $f \rest A = \id A$.
		\item[{\rm(b)}] For every $A \loe G$, $A \in \Ef$, there are an isomorphism $\map h{A'}A$ and $B \goe A'$, $B \in \Ef$, such that for every $X \in \Ef$ with $B \loe X$ there exists an embedding $\map f X G$ extending $h$.
		\item[{\rm(c)}] For every $A \in \Ef$, for every embedding $\map e A G$ there is $B \goe A$, $B \in \Ef$, such that for every $X \in \Ef$ with $B \loe X$ there exists an embedding $\map f X G$ extending $e$.
	\end{enumerate}
\end{prop}

\begin{pf}
	(a)$\implies$(b) Take $h := \id A$.
	
	(b)$\implies$(c) Note that the assertion of (b) holds when the isomorphism $h$ is replaced by $h \cmp g$, where $\map{g}{A''}{A'}$ is an arbitrary isomorphism. This follows immediately from Lemma~\ref{LMliftinisos} applied to $g^{-1}$.
	Thus, one can take $g$ to be any automorphism of $A'$, therefore $h \cmp g$ can be an arbitrary embedding of $A'$ into $G$ whose image is $A$, which shows (c).		
	
	(c)$\implies$(a) Take $e := \id A$ and apply (c) obtaining a suitable $B \goe A$. In particular, there is an embedding $\map g B G$ that is identity on $A$. Then $\img g B \in \Ef$ and $A \loe \img g B \loe G$. Fix $X \in \Ef$ with $\img g B \loe X$.
	By Lemma~\ref{LMliftinisos}, there are $X' \goe B$ and an isomorphism $\map{g'}{X'}{X}$ extending $g$. By (c), there is an embedding $\map{f'}{X'}{G}$ that is identity on $A$.
	Finally, $f := f' \cmp (g')^{-1}$ is an embedding of $X$ into $G$ that is identity on $A$.
\end{pf}

We shall say that $G \in \ciagi \Ef$ is \emph{weakly $\Ef$-injective} if it is $\Ef$-universal and satisfies any of the equivalent conditions of Proposition~\ref{PropWEPs}.
Furthermore, we shall say that $G$ is \emph{weakly injective} if it is weakly $\Ef$-injective with $\Ef = \Age{G}$.

\begin{tw}\label{ThmEwaNot}
	Suppose Eve does not have a winning strategy in $\BMG{\Ef}{G}$.
	Then $G$ is weakly $\Ef$-injective.
\end{tw}

\begin{pf}
	First of all, note that $G$ is $\Ef$-universal, since otherwise there would be $A \in \Ef$ not embeddable into $G$ and Eve would have a winning strategy, simply starting the game with $A_0 := A$.
	Assume, for a contradiction, that $G$ is not weakly $\Ef$-injective, we shall use condition (b) of Proposition~\ref{PropWEPs}. Namely, suppose (b) fails and fix a witness $A \loe G$, $A \in \Ef$. We shall describe a winning strategy for Eve.
	Note that the following condition is fulfilled.
	\begin{enumerate}
		\item[($\times$)] For every isomorphism $\map{h}{A'}{A}$, for every $B \goe A'$ with $B \in \Ef$, there exists $B' \goe B$ with $B' \in \Ef$ such that no embedding of $B'$ into $G$ extends $h$.
	\end{enumerate}
	Eve starts with $A_0 := A$.
	Suppose $A_0 \loe A_1 \loe \dots \loe A_{n-1}$ are initial steps of the game $\BMG{\Ef}{G}$, where $n$ is even.
	Eve chooses an isomorphism $h_n$ whose domain is a substructure of $A_{n-1}$ and whose range is $A$. Then she responds with $A_n := B'$ from condition ($\times$) applied to $h := h_n$ and $B := A_{n-1}$. By this way no embedding of $A_n$ into $G$ extends $h_n$.
	This describes Eve's strategy. Note that at each step there are countably many possibilities for choosing an isomorphism onto $A$, therefore an easy book-keeping makes sure that Eve considers all of them. By this way she wins, as in the end no embedding of $\bigcup_{\ntr}A_n$ into $G$ can contain $A$ in its image.
\end{pf}

We will strengthen the above result in Theorem~\ref{THMsoddaof} below.

\begin{tw}\label{ThmOddWinnings}
	Suppose $G \in \ciagi \Ef$ is weakly $\Ef$-injective. Then Odd has a winning strategy in $\BMG{\Ef}{G}$.
\end{tw}

\begin{pf}
	Let $\ciag v$ enumerate a fixed set of generators of $G$. We shall use condition (c) of Proposition~\ref{PropWEPs}, knowing that $G$ is $\Ef$-universal.
		
	Suppose $A_0 \loe \dots \loe A_{n-1}$ form an initial part of $\BMG{\Ef}{G}$ and $n$ is odd.
	We assume that on the way Odd had considered $A_i' \in \Ef$ such that $A_i \loe A_i' \loe A_{i+1}$ for each even $i<n-2$, and he has recorded embeddings $\map{e_i}{A_i'}{G}$ such that $e_i$ extends $e_{i-2}$ and
	$$\{v_0, \dots, v_i\} \subs \img{e_i}{A_i'},$$
	again for each even $i<n-2$.
	Furthermore, if $n>2$, we assume that for every $X \in \Ef$ with $X \goe A_{n-2}$ there exists an embedding $\map{e}{X}{G}$ extending $e_{n-3}$.	
	We now describe Odd's response.
	
	Namely, Odd first finds a copy $B_{n-1} \loe G$ of $A_{n-1}$ and, using Lemma~\ref{LMliftinisos} together with our inductive assumption, finds $A_{n-1}' \goe A_{n-1}$ so that there is an embedding $\map{e_{n-1}}{A_{n-1}'}{G}$ extending $e_{n-3}$ (unless $n=1$), whose image contains $v_i$ for every $i<n$. Odd responds with $A_n \goe A_{n-1}'$ such that the assertion (c) of Proposition~\ref{PropWEPs} holds with $A := A_{n-1}'$, $e := e_{n-1}$, and $B := A_n$. By this way, for every $X \in \Ef$ with $X \goe A_n$, there is an embedding $\map e X G$ extending $e_{n-1}$.
	
	Using this strategy Odd in particular builds an embedding $\map{e_{\infty}}{A_\infty}{G}$, where $A_\infty = \bigcup_{\ntr} A_n = \bigcup_{n \in 2 \Nat} A'_n$ and $e_\infty \rest A'_n = e_n$ for $n \in 2\Nat$. Its image contains the set of generators $\ciag v$, therefore $e_\infty$ is an isomorphism from $A_\infty$ onto $G$.
\end{pf}

Thus, a weakly $\Ef$-injective model in $\ciagi \Ef$ is unique, up to isomorphism.
In~\cite{KubBM} it has been called \emph{$\Ef$-generic} .

\section{Weak amalgamations}

The following concept was introduced and used by Ivanov~\cite{Ivanov}, later by Kechris and Rosendal~\cite{KecRos}, and recently by Kruckman~\cite{KruckPhD}. Ivanov called it the \emph{almost amalgamation property}.

\begin{df}\label{DfWAPs}
	Let $\Ef$ be a class of finitely generated structures. We say that $\Ef$ has the \emph{weak amalgamation property} (briefly: \emph{WAP}) if for every $Z \in \Ef$ there is $Z' \in \Ef$ containing $Z$ as a substructure and such that for every embeddings $\map f {Z'}X$, $\map g {Z'}Y$ with $X,Y \in \Ef$ there exist embeddings $\map{f'}{X}{V}$, $\map{g'}{Y}{V}$ with $V \in \Ef$, satisfying
	$$f' \cmp f \rest Z = g' \cmp g \rest Z.$$	
	We say that $\Ef$ has the \emph{cofinal amalgamation property} (briefly: \emph{CAP}) if $$f' \cmp f = g' \cmp g$$ holds in the definition above.
	Finally, $\Ef$ has the \emph{amalgamation property} (briefly: \emph{AP}) if $Z' = Z$ in the definition above.	
	A subclass $\Ef'$ of a class $\Ef$ is called \emph{cofinal} if $\Ef \subs \her {\Ef'}$.
\end{df}

One needs to admit that the cofinal amalgamation property (which perhaps belongs to the folklore) had been considered earlier by Truss~\cite{Truss}. Obviously, CAP implies WAP and AP implies CAP. Note also that the CAP is equivalent to the existence of a cofinal subclass with the AP.
Finite cycle-free graphs provide an example of a hereditary class satisfying CAP and not AP.

\begin{prop}\label{PropCzteryDwa}
	Let $\Ef$ be a class of structures such that there exists $G \in \ciagi \Ef$ that is weakly $\Ef$-injective. Then $\Ef$ has the weak amalgamation property.	
\end{prop}

\begin{pf}
	Fix $Z \in \Ef$. We may assume that $Z \loe G$.
	We shall use condition (a) of Proposition~\ref{PropWEPs}.
	Namely, let $Z' \in \Ef$ be such that $Z \loe Z' \loe G$ and (a) holds with $A := Z$, $B := Z'$.
	Fix embeddings $\map{f}{Z'}{X}$, $\map{g}{Z'}{Y}$ with $X, Y \in \Ef$.
	Applying Proposition~\ref{PropWEPs} twice, we obtain embeddings $\map{f'}{X}{G}$, $\map{g'}{Y}{G}$ such that both $f' \cmp f$ and $g' \cmp g$ are identity on $Z$. In particular, $f' \cmp f \rest Z = g' \cmp g \rest Z$.
\end{pf}

\begin{df}
	Let $\Ef$ be a class of countable finitely generated structures. We shall say that $\Ef$ is a \emph{weak \fra\ class} if it has JEP, WAP, and contains countably many isomorphic types.
\end{df}

\begin{wn}\label{WnSlabe}
	Let $G$ be a countable structure, let $\Ef$ be a cofinal subclass of $\Age G$, and assume Eve does not have a winning strategy in $\BMG{\Ef}{G}$. Then $\Ef$ is a weak \fra\ class.
\end{wn}

\begin{pf}
	That $\Ef$ has JEP and contains countably many types is the statement of Proposition~\ref{Pjepctmts}. The WAP follows directly from Theorem~\ref{ThmEwaNot} and Proposition~\ref{PropCzteryDwa}.
\end{pf}

It remains to show that every weak \fra\ class has its limit, that is, a suitable countably generated structure with the weak extension property. According to the remark on page 320 in \cite{KecRos}, this can be ``carried over without difficulty" adapting the \fra\ theory presented in the book of Hodges~\cite{Hodges}. One cannot disagree with such a statement, however we are not aware of any text where it has been done explicitly, therefore we present details in the next section.

One needs to admit that Chapter 4 of the recent Ph.D. thesis of Kruckman~\cite{KruckPhD} studies the concept of \emph{generic limits}, defined in topological terms, which turns out to be equivalent to ours.
Finally, the short note of Pabion~\cite{Pab} deals with pre-homogeneous models (which are generic limits), stating their uniqueness without proof, saying that it is an adaptation of \fra's technique, not mentioning the weak amalgamation property at all.

\section{Limits of weak \fra\ classes}

Let $\Ef$ be as above, $Z \in \Ef$, and let $Z' \in \Ef$ be such that $Z \loe Z'$. We shall say that $Z'$ is $Z$-\emph{good} if it satisfies the assertion of Definition~\ref{DfWAPs}, namely, for every embeddings $\map f {Z'}X$, $\map g {Z'}Y$ with $X,Y \in \Ef$ there exist embeddings $\map{f'}{X}{V}$, $\map{g'}{Y}{V}$ with $V \in \Ef$, satisfying
$f' \cmp f \rest Z = g' \cmp g \rest Z$.
Note that WAP says: for every $Z \in \Ef$ there is $Z' \in \Ef$ such that $Z \loe Z'$ and $Z'$ is $Z$-good, while CAP means: for every $Z \in \Ef$ there is $Z' \in \Ef$ such that $Z \loe Z'$ and $Z'$ is $Z'$-good.
Note also that if $Z \loe Z' \loe Z''$ and $Z'$ is $Z$-good then so is $Z''$.

We are now ready to prove the main result of this section.

\begin{tw}
	Let $\Ef$ be a weak \fra\ class. Then there exists a unique, up to isomorphisms, structure $G \in \ciagi \Ef$ that is weakly $\Ef$-injective, and such that $\Ef$ is cofinal in $\Age G$.
	Conversely, if $G$ is a countable weakly injective structure then every cofinal subclass of $\Age G$ is a weak \fra\ class.
\end{tw}

The structure $G$ from the first statement will be called the \emph{limit} of $\Ef$.

\begin{pf}
	Note that the second (``conversely") part is the combination of Propositions~\ref{Pjepctmts} and~\ref{PropCzteryDwa}. It remains to show the first part. Uniqueness follows from Proposition~\ref{PrpUnqns}, therefore it remains to show the existence.
	The construction will rely on the following very simple fact (sometimes called the \emph{Rasiowa-Sikorski Lemma}), well known in forcing theory.
		
	\begin{claim}\label{Lmerbginwr}
		Given a partially ordered set $\poset = \pair P \loe$, given a countable family $\Dee$ of cofinal subsets of $\poset$, there exists an increasing sequence $\Ef = \setof{p_n}{\Ntr}$
		such that $D \cap \Ef \nnempty$ for every $D \in \Dee$.
	\end{claim}
		
	Recall that $D$ is \emph{cofinal} in $\poset$ if for every $p \in P$ there is $q \in D$ such that $p \loe q$.
		
	Now, fix a weak \fra\ class $\Ef$ of finitely generated structures.
	First, we ``localize'' $\Ef$: Namely, we assume that each $A \in \Ef$ lives in the set $\Nat$ of nonnegative integers and the complement $\Nat \setminus A$ is infinite.
	Define the following poset $\poset$.
	The universe of $\poset$ is the class $\Ef$ (refined as above) while $\loe$ is the usual relation of ``being a substructure".
	Define
	$$E_C = \setof{ X \in \Ef }{ C \text{ embeds into } X },$$
	where $C \in \Ef$.
	By assumption, there are countably many sets of the form $E_C$ and the joint embedding property implies that each $E_C$ is cofinal in $\poset$.
	Next, given $A \loe A'$ in $\Ef$ such that $A'$ is $A$-good, given an embedding $\map f {A'} B$, define
	$$D_{A,f} = \setof{ X \in \Ef }{ \text{If } A' \loe X \text{ then }(\exists\;\text{ an embedding } \map g B X)\;(\forall\; x\in A)\; g(f(x))=x }.$$
	By assumption, there are countably many sets of the form $D_{A,f}$. Each of them is cofinal, because of the weak amalgamation property.
	We only need to remember that all structures in $\Ef$ are co-infinite in $\Nat$, so that we always have enough space to enlarge them.
		
	Finally, a sequence 
	$$U_0 \loe U_1 \loe U_2 \loe U_3 \loe \dots$$
	obtained from Claim~\ref{Lmerbginwr} to our family of cofinal sets produces a structure $U_\infty = \bigcup_{\Ntr}U_n$ in $\ciagi \Ef$ that is weakly $\Ef$-injective. In particular, $\Ef$ is cofinal in $\Age{U_\infty}$.
\end{pf}

One of the most important features of \fra\ limits is \emph{universality}. Namely, if $\Ef$ is a \fra\ class and $G$ is its limit then every $X \in \ciagi \Ef$ embeds into $G$. This is not true for weak \fra\ classes (see Examples~\ref{ExNonUnivWFs} and \ref{EXdkjfnodag} below), however the following weaker statement holds true.

\begin{tw}\label{ThmSEGFbaefqw}
	Let $\Ef$ be a weak \fra\ class and let $G$ be its limit. Then for every chain
	$$X_0 \loe X_1 \loe X_2 \loe \cdots$$
	of structures in $\Ef$ such that $X_{n+1}$ is $X_n$-good for each $\ntr$, the union $\bigcup_{\ntr}X_n$ embeds into $G$.
\end{tw}

\begin{pf}
	Let us play the game $\BMG{\Ef}{G}$ in such a way that Odd uses his winning strategy. We shall describe Eve's strategy leading to an embedding of $X = \bigcup_{\ntr}X_n$ into $G$.
	
	Eve starts with any $A_0 \in \Ef$ for which there is an embedding $\map{e_1}{X_1}{A_0}$. Odd responds with $A_1 \goe A_0$.
	Eve uses the WAP in order to find $A_2 \goe A_1$ and an embedding $\map{e_2}{X_2}{A_2}$ so that $e_2 \rest X_0 = e_1 \rest X_0$.
	In general, when $n=2k$ and the position of the game is $A_0 \loe \dots \loe A_{n-1}$, we assume that Eve has already recorded embeddings $\map{e_i}{X_i}{A_{2i-2}}$ satisfying
	\begin{equation}
		e_{i+1} \rest X_{i-1} = e_i \rest X_{i-1}
		\tag{$\star$}\label{EqStarJeden}
	\end{equation}
	for each $i\loe k$.
	Eve responds with $A_n = A_{2k} \goe A_{n-1}$ using the WAP, so that there is an embedding $\map{e_k}{X_k}{A_{2k}}$ satisfying $e_k \rest X_{k-2} = e_{k-1} \rest X_{k-2}$.
	By this way, after infinitely many steps of the game Eve has recorded embeddings $\map{e_i}{X_i}{A_{2i-2}}$ satisfying (\ref{EqStarJeden}) for every $i \in \nat$.
	Define $e = \bigcup_{i \in \nat} e_i \rest X_{i-1}$.
	Then $\map e X {\bigcup_{\ntr}A_n}$ is a well-defined embedding, because of (\ref{EqStarJeden}) and $\bigcup_{\ntr}A_n \iso G$, because Odd was using his winning strategy.
\end{pf}

Note that the above result gives the well known universality of \fra\ limits. This is because if $\Ef$ is a \fra\ class then every $X \in \Ef$ is $X$-good, therefore Theorem~\ref{ThmSEGFbaefqw} applies to every structure in $\ciagi \Ef$.
Below is the announced example showing that the result above cannot be improved.

\begin{ex}\label{ExNonUnivWFs}
	Let $\Ef$ be the class of all finite linear graphs, that is, finite graphs with no cycles and of vertex degree $\loe2$. Connected graphs form a cofinal subclass with the AP, therefore $\Ef$ is a weak \fra\ class. Its limit is $\Zee$, the integers with the linear graph structure (each integer $n$ is connected precisely to $n+1$ and $n-1$).
	The graph $X$ consisting of two disjoint copies of $\Zee$ is in $\ciagi{\Ef}$, however it cannot be embedded into $\Zee$.
\end{ex}

Another important feature of \fra\ limits is \emph{homogeneity}, namely, every isomorphism between finitely generated substructures extends to an automorphism.
This is obviously false in the case of weak \fra\ classes (see Example~\ref{ExNonUnivWFs} above: homogeneity totally fails for disconnected subgraphs of $\Zee$). On the other hand, the following weaker result is true.

\begin{tw}
	Let $\Ef$ be a weak \fra\ class with limit $G$.
	Then for every $A \loe A' \loe G$ such that $A, A' \in \Ef$ and  $A'$ is $A$-good, for every embedding $\map{e}{A'}{G}$ there exists an automorphism $\map h G G$ such that $h \rest A = e \rest A$.
\end{tw}

\begin{pf}
	The proof is a suitable adaptation of the classical back-and-forth argument.
	Namely, let $A_0 := A$, $A_1 := A'$, $f_1 := e$. Let $B_1 \loe G$ be in $\Ef$, such that $\img {f_1}{A_1} \subs B_1$.
	Choose $B_2 \in \Ef$ such that $B_1 \loe B_2 \loe G$ and $B_2$ is $B_1$-good. Applying the weak extension property, find $\map{g_2}{B_2}{G}$ such that $g_2 \cmp f_1 \rest A_0$ is identity. Let $A_2\in\Ef$ be such that $\img{g_2}{B_2} \subs A_2 \loe G$. Enlarging $A_2$ if necessary, we may assume that it is $A_1$-good.
	Inductively, we construct two chains of structures in $\Ef$
	$$A_0 \loe A_1 \loe A_2 \loe \dots \loe G \oraz B_1 \loe B_2 \loe B_3 \loe \dots \loe G$$
	and embeddings $\map{f_{2n+1}}{A_{2n+1}}{B_{2n+1}}$ and $\map{g_{2n}}{B_{2n}}{A_{2n}}$ satisfying for every $\ntr$ the following conditions:
	\begin{enumerate}[itemsep=0pt]
		\item[(1)] $A_{n+1}$ is $A_n$-good and $B_{n+1}$ is $B_n$-good,
		\item[(2)] $g_{2n} \cmp f_{2n-1}$ is identity on $A_{2n-2}$,
		\item[(3)] $f_{2n+1} \cmp g_{2n}$ is identity on $B_{2n-1}$,
		\item[(4)] $\bigcup_{\ntr}A_n = G = \bigcup_{\ntr}B_{n+1}$.
	\end{enumerate}
	Given $f_{2n-1}$ and $g_{2n-2}$, we find $g_{2n}$ and $f_{2n+1}$ exactly in the same way as in the first step, using the weak extension property of $G$.
	Note that we have a freedom to enlarge $A_{2n+2}$ and $B_{2n+1}$ as much as we wish, therefore we can easily achieve (4), knowing that $G$ is countably generated.	
	Now observe that
	$$f_{2n-1} \rest A_{2n-2} = f_{2n+1} \cmp g_{2n} \cmp f_{2n-1} \rest A_{2n-2} = f_{2n+1} \rest A_{2n-2},$$
	because $\img{f_{2n-1}}{A_{2n-2}} \subs B_{2n-1}$ and hence we were able to apply (3) and (4).
	It follows that
	$$f_\infty := \bigcup_{\ntr} f_{2n+1}\rest A_{2n}$$
	is a well-defined embedding of $G$ into itself.
	The same argument shows that
	$$g_\infty := \bigcup_{\ntr} g_{2n+2} \rest B_{2n+2}$$
	is a well-defined embedding of $G$ into itself.
	Conditions (2), (3) make sure that $f_\infty \cmp g_\infty = \id G = g_\infty \cmp f_\infty$, showing that $f_\infty$ is an isomorphism.
	Finally, $$f_\infty \rest A = f_\infty \rest A_0 = f_1 \rest A_0 = e \rest A.$$
	This completes the proof.
\end{pf}

The property stated in the theorem above is called \emph{prehomogeneity} in Pabion~\cite{Pab}.
Note that, again, if $\Ef$ is a \fra\ class then the result above gives full homogeneity, namely, that every embedding between finitely generated substructures extends to an automorphism.

\section{Universality vs. weak amalgamation property}\label{SeccjaSesta}

Let $\Ef$ and $\ciagi \Ef$ be as before. Let $\BMGE{\Ef}{W}$ be the game $\BM{\Ef}$ (namely, the same rules as before) while now \emph{Odd wins} if the resulting structure embeds into $W$. Clearly, if Odd has a winning strategy in $\BMGE{\Ef}{W}$ then $\Ef$ has the JEP and CMT (countably many isomorphic types).

\begin{tw}\label{THMsoddaof}
	Assume $\Ef$ fails the WAP. Then Eve has a winning strategy in $\BMGE{\Ef}{W}$ for every $W \in \ciagi{\Ef}$.
\end{tw}

\begin{pf}
	Fix $W \in \ciagi{\Ef}$. We may assume that $W$ is $\Ef$-universal, since otherwise Eve may start the game with any $A_0 \in \Ef$ not embeddable into $W$, obviously winning the game.
	
	Fix $A \in \Ef$ and suppose $\Ef$ fails the WAP at $A$. 
	Let $\setof{e_n}{\ntr}$ be a fixed enumeration of all embeddings of $A$ into $W$.
	Given an embedding $\map i A B$ with $B \in \Ef$, define $\al(i)$ to be the minimal $k$ such that $j \cmp i = e_k$ for some embedding $\map j B W$.
	This is well-defined because $W$ is $\Ef$-universal and hence contains at least one isomorphic copy of $B$.
	
	Note that, by the failure of WAP at $A$ and by the $\Ef$-universality of $W$, for every $A' \in \Ef$ with $A \loe A'$ there exist embeddings $\map f {A'} X$, $\map g {A'} Y$ with $X,Y \in \Ef$ and such that $f,g$ cannot be amalgamated over $A$, namely, no embeddings $\map {f'}X W$, $\map {g'}Y W$ satisfy $f' \cmp f \rest A = g' \cmp g \rest A$. In particular, $\al(f\rest A) \ne \al(g \rest A)$. We shall use this observation when designing Eve's winning strategy.
	
	Namely, Eve starts the game with $A_0 := A$.
	Suppose $A_0 \loe \dots \loe A_{n-1}$ have already been chosen by the players and assume $n$ is even.
	By the above observation, there are pairs of embeddings $f,g$ of $A_{n-1}$ into structures of $\Ef$ that cannot be amalgamated over $A$. Let $f_n, g_n$ be such a pair in which $\al(f_n \rest A)$ is minimal possible. In particular, $\al(g_n \rest A) > \al(f_n \rest A)$. Without loss of generality, we may assume that $g_n$ is inclusion $A_{n-1} \subs C$, that is, $A_{n-1} \loe C$.
	Eve's response is $A_n := C$. This finishes the description of Eve's strategy.
	
	In order to complete the proof, it suffices to show that there is no embedding of $A_\infty := \bigcup_{\ntr} A_n$ into $W$.
	Towards a contradiction, suppose $\map h {A_\infty} W$ is an embedding. Then $h \rest A = e_\ell$ for some $\ell \in \nat$.
	Note that $\ell > \al(f_n \rest A)$ for each even $n$, because $h \rest A = (h \rest A_n) \cmp g_n \rest A$.
	Finally, notice that $\sett{\al(f_n)}{n \in 2\Nat}$ is strictly increasing, therefore $\ell$ is not a natural number, a contradiction.
\end{pf}

Given a family $\Yu \subs \ciagi \Ef$, denote by $\sbs \Yu$ the family of all structures embeddable into some $U \in \Yu$. Thus, Odd wins in the game $\BMGE{\Ef}{\Yu}$ if the resulting structure embeds into some $U \in \Yu$.

\begin{tw}\label{THMbdsqgvbgo}
	Let $\Yu \subs \ciagi \Ef$ be countable and assume Odd has a winning strategy in $\BMGE{\Ef}{\Yu}$. Then $\Ef$ has the WAP.
\end{tw}

\begin{pf}
	Suppose $\Ef$ fails the WAP at $A$. Then for each $B \in \Ef$ with $A \subs B$ there exist $X,Y \in \Ef$ such that $B \loe X$, $B \loe Y$ and $X,Y$ cannot be amalgamated over $A$, that is, if $\map f X C$, $\map g Y C$ are embeddings with $C \in \Ef$ then $f \rest A \ne g \rest A$.
	Let $\Sigma$ be Odd's winning strategy in $\BMGE{\Ef}{\Yu}$ and let $2^{<\omega}$ denote the tree of all finite binary sequences.
	We can easily construct families $\sett{A_s}{s\in 2^{<\omega}} \subs \Ef$ and $\sett{B_s}{s\in 2^{<\omega}} \subs \Ef$ satisfying
	\begin{enumerate}[itemsep=0pt]
		\item[(0)] $A_\emptyset = A$.
		\item[(1)] $A_s \loe B_s \loe A_{s\concat i}$ for $i<2$.
		\item[(2)] $A_{s\concat0}$, $A_{s\concat1}$ cannot be amalgamated over $A$.
		\item[(3)] $B_s$ is Odd's response, according to $\Sigma$, to the sequence
		$$A_\emptyset \loe B_\emptyset \loe A_{s \rest 1} \loe B_{s \rest 1} \loe \dots \loe A_s.$$
	\end{enumerate}
	This is done by an easy induction: Having defined $A_s$, clause (3) provides the definition of $B_s$; given $B_s$ it is always possible to find $A_{s\concat0}$, $A_{s\concat1}$ satisfying (2).
	We may additionally assume $B_s = A_{s\concat0} \cap A_{s\concat1}$.
	
	Now, given $\sig \in 2^\omega$, let $A_\sig = \bigcup_{\ntr} A_{\sig \rest n}$.
	Then $A_\sig \in \sig \Ef$ and since it is a result of a play in which Odd was using his winning strategy $\Sigma$, we conclude that $A_\sig$ embeds into some $U_\sig \in \Yu$.
	By assumption, $\Yu$ is countable, therefore there is $U \in \Yu$ such that the set
	$$S = \setof{\sig \in 2^\omega}{U_\sig = U}$$ is uncountable.
	For each $\sig \in S$ choose an embedding $\map{e_\sig}{A_\sig}{U}$.
	Since there are countably many embeddings of $A$ into $U$ (recall that $U \in \ciagi{\Ef}$ is countable), there are $\sig \ne \tau$ in $S$ such that $e_\sig \rest A = e_\tau \rest A$.
	Let $s = \sig \cap \tau$, namely, $s \in 2^{<\omega}$ is maximal below $\sig$ and $\tau$. We may assume $s\concat0 \subs \sig$ and $s\concat1 \subs \tau$. Then $e_\sig \rest A_{s\concat0}$ and $e_\tau \rest A_{s\concat1}$ are embeddings violating (2). This is a contradiction.
\end{pf}

Note that one can relax the assumption on $\Yu$ in Theorem~\ref{THMbdsqgvbgo} by requiring that $\Yu$ is of cardinality strictly less than the continuum. In particular, we obtain the following

\begin{wn}
	Assume there exists a family $\Yu \subs \ciagi \Ef$ of cardinality $< 2^{\aleph_0}$ such that each $X \in \ciagi{\Ef}$ embeds into some $U \in \Yu$.
	Then $\Ef$ has the WAP.
\end{wn}

\begin{wn}
	If $\ciagi{\Ef}$ has a universal model then $\Ef$ is a weak \fra\ class.
\end{wn}

Finally, note that a $\ciagi{\Ef}$-universal model may not be weakly $\Ef$-injective and a weakly $\Ef$-injective model may not be $\ciagi{\Ef}$-universal.

\begin{ex}[{cf. \cite[Ex. 11]{KubBM}}]\label{EXdkjfnodag}
	Fix $k>1$ and let $\Gee_k$ be the class of all finite graphs with vertex degree $\loe k$.
	Obviously, $\Gee_k$ has the JEP and countably many types.
	Recall that a graph is $k$-regular if each vertex has degree precisely $k$. A standard argument shows that every graph in $\Gee_k$ embeds into a finite $k$-regular graph. Furthermore, if $G \in \ciagi \Gee_k$ is connected, $k$-regular, and $G$ is a subgraph of $H \in \ciagi \Gee_k$ then $G$ is a component of $H$.
	It follows that $\Gee_k$ has the cofinal amalgamation property. Indeed, if $Z \in \Gee_k$ is $k$-regular and $\map f Z X$, $\map g Z Y$ are embeddings with $X, Y \in \Gee_k$ then the free amalgamation of $f$ and $g$ is a graph of vertex degree $\loe k$, therefore a member of $\Gee_k$.
	Let $U_k$ denote the countable weakly $\Gee_k$-injective graph (the generic limit of $\Gee_k$).
	Then $U_k$ can be described by the following properties:
	\begin{enumerate}[itemsep=0pt]
		\item[(1)] Each component of $U_k$ is finite and $k$-regular.
		\item[(2)] For every finite $k$-regular graph $G$ infinitely many components of $U_k$ are isomorphic to $G$.
	\end{enumerate}
	Indeed, Odd has a simple winning strategy in $\BMG{\Gee_k}{U_k}$. Namely, at stage $n$ he first enlarges each component to a (finite) $k$-regular graph and next, by adding new components, he makes sure that the first $n$ isomorphic types of finite $k$-regular graphs appear, each of them at least $n$ times. Of course, the strategy depends on a fixed enumeration of all finite $k$-regular graphs.
	
	Now let $k=2$. Note that $U_2$ is the disjoint union of finite cycles and for every $n>2$ there are infinitely many cycles of length $n$.
	On the other hand, $\ciagi \Gee_2$ contains infinite $2$-regular graphs, for example, the integers $\Zee$ in which $a,b$ are connected if and only if $|a-b|=1$.
	Note also that $\ciagi \Gee_2$ has a unique, up to isomorphism, universal regular graph: the disjoint sum of $U_2$ and infinitely many isomorphic copies of $\Zee$.
	
	Now let $k=3$. We claim that there are continuum many pairwise non-isomorphic connected countable $3$-regular graphs. Consequently, $\ciagi \Gee_3$ 
	has no universal graph. Even more: the minimal cardinality of a family $\Yu \subs \ciagi \Gee_3$ such that every graph in $\Gee_3$ embeds into some graph of $\Yu$ has cardinality continuum.
	This is because, each connected $3$-regular graph embeds onto a component of any graph in $\ciagi \Gee_3$ and therefore each of the continuum many pairwise non-isomorphic components must appear in the graphs from $\Yu$.
	
	It remains to construct the announced $3$-regular graphs.
	Fix $A \subs \Nat$ with $\min A > 3$. Let $T_A$ be a countable infinite tree (a connected cycle-free graph) such that
	$$A = \setof{\deg(v)}{v \in T}.$$
	Now replace each vertex $v \in T$ by a cycle $C_v$ of length $\deg(v)$ and modify the graph so that it becomes $3$-regular.
	Denote the resulting graph by $G_A$.
	Obviously, if $A \ne B$ then $G_A$ is not isomorphic to $G_B$.
\end{ex}

The last example shows that there may be no universal model in $\ciagi \Ef$ even if $\Ef$ is a weak \fra\ class with the cofinal AP.
The results above show that it is not so obvious to find classes of finite models with no WAP. Namely, if $\Ef$ has the JEP and CMT then the failure of WAP implies, in particular, no universal model in $\ciagi \Ef$. Nevertheless, below is an example of a class of finite graphs with JEP that fails the WAP.

\begin{ex}
	Let $\Gee$ be the class of all finite graphs in which different cycles of the same length are disjoint. Recall that a \emph{cycle} is a path $\{v_0,\dots, v_n\}$ such that $v_0 = v_n$ and $v_i \ne v_j$ whenever $0 \loe i < j < n$.
	Clearly, $\Gee$ has the JEP, as it is closed under disjoint sums. We claim that $\Gee$ fails the WAP.
	
	Fix a graph $G \in \Gee$ containing three different vertices $a,b,c$ such that $a$ is adjacent to $b$ and $c$.
	Suppose $G' \in \Gee$ is such that $G \loe G'$.
	Fix an integer $k$ greater than the cardinality of $G'$.
	Let $X$ be an extension of $G'$ by adding a path of length $2k$ starting at $a$ and finishing at $b$. Then $X \in \Gee$, because $G'$ has no cycle of length $\goe k$.
	Let $Y$ be a similar extension of $G'$ but now we add one more edge, namely, we connect $c$ with  the unique vertex $y$ on the new path whose distance to $a$ is $k-1$.
	Again, $Y \in \Gee$, by the same reason as before.
	Suppose $\map fXW$, $\map gYW$ are embeddings such that $f \rest G = g \rest G$.
	The cycles $C_X = \dn ab \cup (X \setminus G')$ and $C_Y = \dn ab \cup (Y \setminus G')$ have both length $2k+1$. Their images $\img f{C_X}$, $\img g{C_Y}$ have a common edge $\dn ab$ but they cannot be equal because then $g(y) = f(x)$, where $x$ is the unique vertex in $X \setminus G'$ whose distance to $a$ is $k-1$. On the other hand, $g(y)$ is adjacent to $f(c) = g(c)$ and $f(x)$ is not adjacent to $f(c) = g(c)$, because $f$, $g$ are embeddings.
	It follows that $W \notin \Gee$.
\end{ex}

\section{Final remarks}

We can now summarize results involving the abstract Banach-Mazur game.

\begin{wn}\label{CORwnergbreg}
	Let $\Ef$ be a class of finitely generated structures with the joint embedding property and with countably many isomorphic types.
	\begin{enumerate}[itemsep=0pt]
		\item[{\rm(1)}] If there is a family $\Yu \subs \ciagi \Ef$ of cardinality strictly less than the continuum, such that Odd has a winning strategy in $\BMGE{\Ef}{\Yu}$, then $\Ef$ has the weak amalgamation property and consequently there exists a unique $U \in \ciagi \Ef$ such that Odd has a winning strategy in the game $\BMG{\Ef}{U}$.
		The structure $U$ is characterized by the weak $\Ef$-injectivity or, equivalently, by the $\Ef$-universality and prehomogeneity.
		\item[{\rm(2)}] If $\Ef$ fails the weak amalgamation property, then Eve has a winning strategy in $\BMGE{\Ef}{W}$ for any $W \in \ciagi \Ef$; her strategy is also winning in the game $\BMG{\Ef}{W}$.
	\end{enumerate}
\end{wn}

Thus, our results show that the weak amalgamation property plays the crucial role in the game.
It is natural to ask the following:

\begin{question}\label{PbRGdg}
	(a) Does there exist a hereditary weak \fra\ class without the cofinal amalgamation property? (b) How about weak \fra\ classes of finite graphs?
\end{question}

Neither Ivanov~\cite{Ivanov} nor Kechris and Rosendal~\cite{KecRos} addressed this question explicitly, all of their examples had the CAP. Nevertheless, the answer to Question~\ref{PbRGdg}(a) is, as one can expect, positive.
It seems that the first (and perhaps the simplest) example is due to Maurice Pouzet, presented in Pabion~\cite{Pab} in 1972. The example reads as follows. Take the set $\Qyu$ of rational numbers and define a ternary relation $R$ by setting $R(x,y,z)$ if and only if $x<y$, $x<z$, and $y \ne z$.
Let $\Ef = \Age{\Qyu, R}$. Then $\Ef$ is a weak \fra\ class with no cofinal amalgamation property. The failure of CAP comes from the fact that on a finite set the relation $R$ is not able to distinguish the ordering between the two greatest elements.
Another example was found by Alex Kruckman in 2015 (email communication), see~\cite[Example 3.4.7 on p. 73]{KruckPhD}. Formally, this example is not hereditary, however its hereditary closure still fails the CAP.
Recently, Aristotelis Panagiotopoulos (email communication) has found another example which lead us to discovering an example answering Question~\ref{PbRGdg}(b).
All these examples will be presented in the subsequent paper~\cite{KKKP}.

Corollary~\ref{CORwnergbreg}(1) suggests that there is a link between WAP and the universality number. Specifically, given a class $\Ef$ of finitely generated structures, denote by $\covsig(\Ef)$ the \emph{universality number} of $\ciagi \Ef$, that is, the minimal cardinality of a family $\Yu \subs \ciagi \Ef$ such that every $X \in \ciagi \Ef$ embeds into some $U \in \Yu$.
Clearly, $\covsig(\Ef)=1$ whenever $\Ef$ is a \fra\ class.
Example~\ref{EXdkjfnodag} shows that $\covsig(\Ef)$ might be equal to the continuum if $\Ef$ is just a cofinal \fra\ class. On the other hand, $\Ef$ has the weak amalgamation property, whenever consistently $\covsig(\Ef) < 2^{\aleph_0}$.
Note that the WAP is absolute (i.e., it does not depend on the universe of set theory in which $\Ef$ lives).

\begin{question}
	Does there exist a hereditary class $\Ef$ of finitely generated structures in a fixed countable signature, with the joint embedding property and with countably many isomorphic types, such that consistently $\aleph_0 < \covsig(\Ef) < 2^{\aleph_0}$?
\end{question}

Such a class necessarily has the WAP and hence it is a weak \fra\ class.
The assumption on countably many types is relevant here, as the following fact shows.

\begin{prop}
	Let $\Ef$ be a class of finitely generated structures of a fixed countable signature. Assme $\Ef$ has the joint embedding property, the amalgamation property, and precisely $\kappa$ many isomorphic types, where $\kappa > \aleph_0$.
	Then
	$$\cov_\sig(\Ef) = \cf \left(\ctbls{\kappa}, \subs \right).$$
\end{prop}

\begin{pf}
	Note that all the structures in $\Ef$ are countable (possibly finite).
	For the purpose of this proof, we refine $\Ef$ by selecting a single copy of each isomorphic type. By this way $\Ef$ becomes a set of cardinality $\kappa$ and moreover isomorphic models in $\Ef$ are equal.
	Let $\lam = \cf(\ctbls{\kappa}, \subs)$.
	We can find countable collections $\Ef_i \subs \Ef$, $i < \lam$, so that every countable subcollection of $\Ef$ is contained in some $\Ef_i$. By a simple closing-off argument, we may assume each $\Ef_i$ has both the joint embedding the amalgamation property.
	Now each $\Ef_i$ is a cofinal \fra\ class, as it has countably many models.
	Let $U_i$ be the \fra\ limit of $\Ef_i$ and let $\Yu = \setof{U_i}{i<\lam} \subs \sig \Ef$.
	Fix $M \in \sig \Ef$ and let $\Gee$ be the collection of all $A \in \Ef$ embeddable into $M$. Find $j<\lam$ such that $\Gee \subs \Ef_j$. By the universality of \fra\ limits, we conclude that $M$ is embeddable into $U_j$.
	This shows that $\cov_\sig(\Ef) \loe \lam$.
	
	In order to show the reversed inequality, fix a family $\Yu \subs \sig \Ef$ such that every structure in $\sig \Ef$ embeds into some $U \in \Yu$.
	For each $U \in \Yu$ define
	$\Gee_U = \setof{A \in \Ef}{A \text{ embeds into }U}.$
	Then $\Gee_U$ is countable (recall that isomorphic structures in $\Ef$ are equal).
	Fix a countable $\Aaa \subs \Ef$. Using the joint embedding property of $\Ef$, find $M \in \sig \Ef$ such that each $A \in \Aaa$ embeds into $M$.
	Next, find $U \in \Yu$ such that $M$ embeds into $U$. In particular $\Aaa \subs \Gee_U$. This shows that $\sett{\Gee_U}{U \in \Yu}$ is cofinal in $(\ctbls{\Ef}, \subs)$.
	As $|\Ef| = \kappa$, this shows that $|\Yu| \goe \lam$.
\end{pf}

Note that the last part of the proof above does not use the amalgamation property.
Note also that $\cf(\ctbls{\aleph_n}, \subs) = \aleph_n$ for every integer $n>0$, therefore $\covsig(\Ef)$ could be any cardinal number strictly between $\aleph_0$ and $\aleph_\omega$.
A concrete example could be the class $\Emm_S$ of all finite metric spaces with distances in a fixed infinite additive subgroup $S \subs \Err$. A metric can be defined by using countably many relations (for each rational number $r$ we have a binary relation saying that the distance is $<r$). Once $S$ is closed under addition, $\Emm_S$ has the amalgamation property. Obviously, it has the JEP and the number of isomorphic types in $\Emm_S$ is precisely 
the cardinality of $S$.
Thus, if $|S| = \aleph_n$ with $\ntr$ then $\covsig(\Emm_S) = \aleph_n$.
In particular, if continuum is large enough, there exist classes $\Ef$ with $\covsig(\Ef)$ equal to any prescribed cardinal $\aleph_n$ below the continuum.

\end{document}